# A summation by Genčev

Donal F. Connon

dconnon@btopenworld.com

5 May 2008

**Abstract**

At a conference at the University of Ostrava in September 2007, Genčev reported that

$$\sum_{n=1}^{\infty} \frac{H_n^{(1)}}{2n+1} \binom{2n}{n}^{-1} = -\frac{\pi}{3\sqrt{3}} \log 3 - \frac{2}{\sqrt{3}} \delta_1$$

where $\delta_1$ is defined by

$$\delta_1 = i\left[ \operatorname{dilog}\left(\frac{\sqrt{3}-i}{2\sqrt{3}}\right) - \operatorname{dilog}\left(\frac{\sqrt{3}+i}{2\sqrt{3}}\right) \right]$$

$$\operatorname{dilog}(x) = Li_2(1-x) = \sum_{n=1}^{\infty} \frac{(1-x)^n}{n^2}$$

In the above $H_n^{(1)}$ is the harmonic number $H_n^{(1)} = \sum_{k=1}^{n} \frac{1}{k}$ and $i = \sqrt{-1}$.

This note indicates a possible approach to the proof.

**Proof**

We have the Beta function $B(u,v)$ defined for $\operatorname{Re}(u) > 0$ and $\operatorname{Re}(v) > 0$ by the Eulerian integral

$$B(u,v) = \int_0^1 t^{u-1}(1-t)^{v-1}\,dt$$

and it is well known that

$$B(u,v) = \frac{\Gamma(u)\Gamma(v)}{\Gamma(u+v)}$$

where $\Gamma(u)$ is the gamma function. We then obtain

$$\int_0^1 (1-v)^n v^n \, dv = \frac{\Gamma(n+1)\Gamma(n+1)}{\Gamma(2n+2)}$$

which gives us

(1) $$\frac{1}{2n+1}\binom{2n}{n}^{-1} = \int_0^1 (1-v)^n v^n \, dv$$

We have the familiar integral for the digamma function [2, p.15]

$$\psi(u+1) + \gamma = \int_0^1 \frac{1-t^u}{1-t} \, dt$$

where integration by parts gives us

$$= (1-t^u)\log(1-t)\Big|_0^1 - u\int_0^1 t^{u-1} \log(1-t) \, dt$$

and thus we see that

$$\psi(u+1) + \gamma = -u \int_0^1 t^{u-1} \log(1-t) \, dt$$

In the case where $u$ is an integer we have the known result

(2) $$H_n^{(1)} = \gamma + \psi(n+1) = -n \int_0^1 t^{n-1} \log(1-t) \, dt$$

Combining (1) and (2) and making the summation gives us

$$\sum_{n=1}^\infty \frac{H_n^{(1)}}{2n+1}\binom{2n}{n}^{-1} = -\sum_{n=1}^\infty n \int_0^1 (1-u)^n \frac{\log u}{1-u} \, du \int_0^1 v^n (1-v)^n \, dv$$

We have the geometric series

$$\sum_{n=1}^\infty w^n = \frac{w}{1-w}$$

and differentiation results in

$$\sum_{n=1}^{\infty} n w^n = w \frac{d}{dw} \frac{w}{1-w} = \frac{w}{(1-w)^2}$$

and we end up with the double integral representation

$$\sum_{n=1}^{\infty} \frac{H_n^{(1)}}{2n+1} \binom{2n}{n}^{-1} = -\int_0^1 \int_0^1 \frac{v(1-v)\log u}{[1-(1-u)v(1-v)]^2} \, du \, dv$$

It is easily seen that

$$\int \frac{\log u}{[1-a(1-u)]^2} \, du = -\frac{[1-a(1-u)]\log[1-a(1-u)] - au \log u}{a(1-a)[1-a(1-u)]}$$

and we obtain the definite integral

$$\int_0^1 \frac{\log u}{[1-a(1-u)]^2} \, du = \frac{\log(1-a)}{a(1-a)}$$

or equivalently

$$\int_0^1 \frac{v(1-v)\log u}{[1-(1-u)v(1-v)]^2} \, du = \frac{\log[1-v(1-v)]}{1-v(1-v)}$$

The problem is then reduced to the single integral

(3) $$\sum_{n=1}^{\infty} \frac{H_n^{(1)}}{2n+1} \binom{2n}{n}^{-1} = -\int_0^1 \frac{\log[1-v(1-v)]}{1-v(1-v)} \, dv$$

The Wolfram Integrator gives us in a scintilla temporis

$$6 \int \frac{\log[1-v(1-v)]}{1-v(1-v)} \, dv = -i\sqrt{3} \log^2\left[v - \sqrt[3]{-1}\right] + \frac{3(-1)^{2/3} \log^2\left[v + (-1)^{2/3}\right]}{1 + \sqrt[3]{-1}}$$

$$-4\sqrt{3} \tan^{-1}\left(\frac{2v-1}{\sqrt{3}}\right)\left[\log\left[v - \sqrt[3]{-1}\right] + \log\left[v + (-1)^{2/3}\right] - \log(v^2 - v + 1)\right]$$

$$-\frac{6(-1)^{2/3}}{1+\sqrt[3]{-1}}\left[\log\left[v+(-1)^{2/3}\right]\log\left(\frac{1+(-1)^{2/3}v}{1+\sqrt[3]{-1}}\right) + Li_2\left(-i\frac{v+(-1)^{2/3}}{\sqrt{3}}\right)\right]$$

$$+\frac{6(-1)^{2/3}}{1+\sqrt[3]{-1}}\left[\log\left[v-\sqrt[3]{-1}\right]\log\left(-(-1)^{2/3}\frac{1+(-1)^{2/3}v}{1+\sqrt[3]{-1}}\right)+Li_2\left(\frac{v+(-1)^{2/3}}{1+\sqrt[3]{-1}}\right)\right]$$

At $v=1$ we find for the right-hand side

$$=-i\sqrt{3}\log^2\left[1-\sqrt[3]{-1}\right]+\frac{3\log^2 2}{1+\sqrt[3]{-1}}$$

$$-4\sqrt{3}\tan^{-1}\left(\frac{1}{\sqrt{3}}\right)\left[\log\left[1-\sqrt[3]{-1}\right]+\log 2\right]$$

$$-\frac{6}{1+\sqrt[3]{-1}}\left[\log 2\log\left(\frac{2}{1+\sqrt[3]{-1}}\right)+Li_2\left(-i\frac{2}{\sqrt{3}}\right)\right]$$

$$+\frac{6}{1+\sqrt[3]{-1}}\left[\log\left[1-\sqrt[3]{-1}\right]\log\left(-\frac{2}{1+\sqrt[3]{-1}}\right)+Li_2\left(\frac{2}{1+\sqrt[3]{-1}}\right)\right]$$

where we have assumed that, of the three possible cubic roots, we are to use $(-1)^{2/3}=1$.

We have

$$\log\left[1+e^{ix}\right]=\log\left[2\cos(x/2)\right]+ix/2$$

$$\sqrt[3]{-1}=\left[\cos\pi+i\sin\pi\right]^{1/3}=\cos\left(\frac{\pi}{3}\right)+i\sin\left(\frac{\pi}{3}\right)$$

and therefore we have

$$\log\left[1+\sqrt[3]{-1}\right]=\log\left[2\cos\left(\frac{\pi}{6}\right)\right]+i\frac{\pi}{6}$$

and $\quad \cos\left(\frac{\pi}{6}\right)=\frac{\sqrt{3}}{2}$

At $v=0$ we find for the right-hand side

$$=-i\sqrt{3}\log^2\left[-\sqrt[3]{-1}\right]-4\sqrt{3}\tan^{-1}\left(\frac{-1}{\sqrt{3}}\right)\left[\log\left[-\sqrt[3]{-1}\right]\right]-\frac{6}{1+\sqrt[3]{-1}}\left[Li_2\left(-i\frac{1}{\sqrt{3}}\right)\right]$$

$$+\frac{6}{1+\sqrt[3]{-1}}\left[\log\left[-\sqrt[3]{-1}\right]\log\left(-\frac{1}{1+\sqrt[3]{-1}}\right)+Li_2\left(\frac{1}{1+\sqrt[3]{-1}}\right)\right]$$

In a personal communication, Jonathan Sondow has kindly informed me that Mathematica produces the approximate value of 0.234163… for both sides of equation (3).

Donal F. Connon
Elmhurst
Dundle Road
Matfield
Kent TN12 7HD